\documentclass{ifacconf}

\usepackage{graphicx}      
\usepackage{natbib}        

\usepackage{physics}
\usepackage{tikz}
\usetikzlibrary{patterns}
\usepackage{subcaption}
\usepackage[version=4]{mhchem}
\newtheorem{remark}{Remark}

\begin{document}
\begin{frontmatter}

\title{A First-Engineering Principles Model for Dynamical Simulation of a Calciner in Cement Production}

\thanks[footnoteinfo]{This work was supported by Innovation Fund Denmark, Ref. 2053-00012B}

\author[DTUComp,FLS]{Jan Lorenz Svensen} 
\author[FLS]{Wilson Ricardo Leal da Silva} 
\author[DTUComp]{John Bagterp J\o rgensen}

\address[DTUComp]{Department of Applied Mathematics and Computer Science, Technical University of Denmark, 2800 Lyngby, Denmark (e-mail: jlsv@dtu.dk, jbjo@dtu.dk).}
\address[FLS]{FLSmidth A/S, 2500 Valby, Denmark (e-mail: wld@flsmidth.com)}

\begin{abstract}                
We present an index-1 differential-algebraic equation (DAE) model for dynamic simulation of a calciner in the pyro-section of a cement plant. The model is based on first engineering principles and integrates reactor geometry, thermo-physical properties, transport phenomena, stoichiometry and kinetics, mass and energy balances, and algebraic volume and internal energy equations in a systematic manner. The model can be used for dynamic simulation of the calciner. We also provide simulation results that are qualitatively correct. The calciner model is part of an overall model for dynamical simulation of the pyro-section in a cement plant. This model can be used in design of control and optimization systems to improve the energy efficiency and \ce{CO2} emission from cement plants.


\end{abstract}


\end{frontmatter}

\section{Introduction}
\label{sec:Introduction}
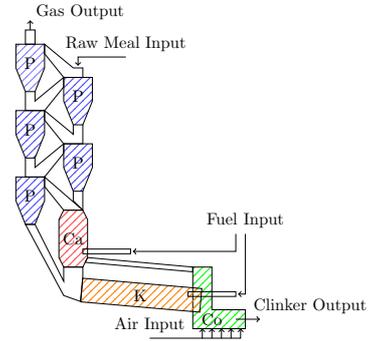
\begin{figure}[tb]
    \centering
    \resizebox{!}{4.5cm}{%
        \begin{tikzpicture}    
            \draw[pattern=north east lines, pattern color=blue!50] (0.2,5.5) -- (0.2,6) -- (0.8,6) -- (0.8,5.5) -- (0.6,5) -- (0.4,5)-- cycle;
            \draw[pattern=north east lines, pattern color=blue!70] (1.2,4.8) -- (1.2,5.3) -- (1.8,5.3) -- (1.8,4.8) -- (1.6,4.3) -- (1.4,4.3)-- cycle;
            \draw[pattern=north east lines, pattern color=blue!70] (0.2,4.1) -- (0.2,4.6) -- (0.8,4.6) -- (0.8,4.1) -- (0.6,3.6) -- (0.4,3.6)-- cycle;
            \draw[pattern=north east lines, pattern color=blue!70] (1.2,3.4) -- (1.2,3.9) -- (1.8,3.9) -- (1.8,3.4) -- (1.6,2.9) -- (1.4,2.9)-- cycle;
            \draw[pattern=north east lines, pattern color=blue!70] (0.2,2.7) -- (0.2,3.2) -- (0.8,3.2) -- (0.8,2.7) -- (0.6,2.2) -- (0.4,2.2)-- cycle;

            \draw  (1.4,5.5) -- (0.8,6) -- (0.8,5.8) -- (1.2,5.3)-- (1.6,5.3) -- (1.6,5.5) -- cycle;
            \draw (0.6,5) -- (0.6,4.8) -- (1.2,5.3) -- (1.2,5.1) -- (0.8,4.6)-- (0.4,4.6) -- (0.4,5);
            \draw (1.4,4.3) -- (1.4,4.1) -- (0.8,4.6) -- (0.8,4.4) -- (1.2,3.9)-- (1.6,3.9) -- (1.6,4.3);
            \draw (0.6,3.6) -- (0.6,3.4) -- (1.2,3.9) -- (1.2,3.7) -- (0.8,3.2)-- (0.4,3.2) -- (0.4,3.6);
            \draw (1.6,2.5)  -- (0.8,3.2) -- (0.8,3.0) -- (1.2,2.5)-- (1.6,2.5);
            \draw (1.4,2.9) -- (1.4,2.7) -- (1.6,2.5) -- (1.6,2.9);
            \draw (0.5,2.8) node{P};
            \draw (1.5,3.5) node{P};
            \draw (0.5,4.2) node{P};
            \draw (1.5,4.9) node{P};
            \draw (0.5,5.6) node{P};
            
            \draw[pattern=north east lines, pattern color=red!70] (1.1,1.5) -- (1.1,2.3) -- (1.2,2.5) -- (1.6,2.5) -- (1.7,2.3) -- (1.7,1.5) -- (1.6,1.3) -- (1.2,1.3)-- cycle;
            \draw (1.4,1.9) node{Ca};
            
            \draw[pattern=north east lines, pattern color=green] (3.9,1+0.3) -- (3.9,1-1) -- (5,1-1) -- (5,1-0.6) -- (4.3,1-0.6) -- (4.3,1+0.3)-- cycle;
            \draw[rotate around={-5:(0,0)},pattern=north east lines, pattern color=orange] (1.5,1-0.3) -- (4,1-0.3) -- (4,1+0.2) -- (1.5,1+0.2) -- cycle;
            \draw (2.8,0.7) node{K};
            
            \draw (3.9,1+0.3) -- (3.9,1+0.2) -- (1.65,1.4) -- (1.7,1.5) -- cycle;
            \draw (1.6,1.3) -- (1.6,1.1) -- (1.55,0.575) -- (1.2,0.7) -- (0.4,2.2) -- (0.6,2.2) -- (1.2, 1) -- (1.2,1.3) -- cycle;
            \draw (4.3,0.2) node{Co};

            \draw  (3.8,0.68) rectangle (4.8,0.78);
            \draw  (1.6,1.58) rectangle (2.6,1.68);
            \draw  (0.4,6) rectangle (0.6,6.3);

            \draw[->] (4.8,0.2) -- node[anchor=south west]{Clinker Output} (5.3,0.2);
            \draw[->] (5.0,2) node[anchor=south]{Fuel Input} -- (5.0,0.73)  -- (4.85,0.73);
            \draw[->] (4.8,2) -- (4.8,1.63)  -- (2.65,1.63);
            \draw[->] (0.5,6.3) -- node[anchor=south west]{Gas Output} (0.5,6.5);
            \draw[->] (2.5,5.73) node[anchor=south]{Raw Meal Input} -- (1.5,5.73)  -- (1.5,5.55);
            \draw[->] (3.0,-0.2) node[anchor=south]{Air Input} -- (4.1,-0.2)  -- (4.1,0.0);
            \draw[->] (4.1,-0.2) -- (4.3,-0.2)  -- (4.3,0.0);
            \draw[->] (4.3,-0.2) -- (4.5,-0.2)  -- (4.5,0.0);
            \draw[->] (4.5,-0.2) -- (4.7,-0.2)  -- (4.7,0.0);
            \draw[->] (4.7,-0.2) -- (4.9,-0.2)  -- (4.9,0.0);
    \end{tikzpicture}    
    }
    \caption{The pyro-section for production of cement clinker in a cement plant consists of preheating cyclones ({\color{blue}P}), a calciner ({\color{red}Ca}), a rotary kiln ({\color{orange}K}), and a cooler ({\color{green}Co}).}\label{fig:production}
\end{figure}

The production of cement clinker is the main source of \ce{CO2} emissions in cement manufacturing. Cement manufacturing is responsible for 8\% of the global \ce{CO2} emissions and about 25\% of all industrial \ce{CO2} emissions \citep{CO2Techreport}. Along with process modifications for carbon capture and storage (CCS), digitalization, control, and optimization are important tools in the transition to zero \ce{CO2}-emission cement plants. Development of such digitalization, control and optimization tools require dynamic simulation and digital twins for the cement plant, and the pyro-section in particular. Mathematical models for dynamic simulation of the pyro-section in cement plants are not available. Fig. \ref{fig:production} illustrates the pyro-section of a cement plant. The pyro-section consists of pre-heating cyclones, a calciner, a rotary kiln, and a cooler.  In this paper, we provide a mathematical model for dynamic simulation of the calciner. A related paper provides a mathematical model for dynamic simulation of the rotary kiln \citep{Svensen2024Kiln}, while papers for the pre-heating cyclones and the cooler are being prepared. Accordingly, the contribution of this paper is a dynamic simulation model for a subunit in the pyro-processing section of a cement plant, namely the calciner. This model is relevant for traditional cement plants as well as modern cement plants designed for carbon capture (oxy-combustion with carbon capture or post carbon capture) and useful for design of control and optimization systems for such plants.

\cite{MUJUMDAR2007} modeled a cyclone-based calciner using a quasi steady state approximation for the coal particles and a dynamic description for the raw meal and gases. \cite{Kahawalage2017} used a CFD approach to model an entrainment calciner. 
\cite{ILIUTA2002805} suggested a 1D dynamic Eulerian model based with detailed combustion kinetics but without a kinetic calcination model. Furthermore, \cite{Kahawalage2017} and \cite{ILIUTA2002805} assume constant heat capacities. Compared to the existing literature, we provide a mathematical 1D model for dynamic simulation of a single elongated chamber calciner (different from a cyclone calciner) that is based on rigorous thermo-physical properties and kinetic expressions for the calcination as well as the combustion. The model is the result of a novel systematic modeling methodology that integrates thermo-physical properties, transport phenomena, and stoichiometry and kinetics with mass and energy balances. The resulting model is a system of index-1 differential algebraic equations (DAEs). 

This paper is organized as follows. Section \ref{sec:Calciner} presents the calciner, while Section \ref{sec:CalcinerModel} describes the mathematical model for the calciner. Section \ref{sec:SimulationResults} presents simulation results and conclusions are provided in Section \ref{sec:Conclusion}.

\section{The calciner}
\label{sec:Calciner}
The main \ce{CO2} contributing reactions in cement manufacturing are calcination of limestone and combustion of coal to provide the heat for the calcination
\begin{subequations}
\begin{alignat}{3}
&\ce{CaCO3} + \ce{heat} \rightarrow \ce{CaO} + \ce{CO2}, \,\, &&\Delta H_r = 179.4\frac{\text{kJ}}{\text{mol}}, \\
&\ce{C} + \ce{O2} \rightarrow \ce{CO2} + \ce{heat}, \, && \Delta H_r = -395.5\frac{\text{kJ}}{\text{mol}}.
\end{alignat}
\end{subequations}
Consequently, the \ce{CO2} emission for the calcination alone excluding other reactions and heating is 0.356 kg \ce{CO2} per produced kg \ce{CaO}. The pyro-section of the cement plant is designed for transfer of heat at high temperatures to the cement raw meal to facilitate the calcination and other cement clinker forming reactions. The cement raw meal has a well controlled chemical composition and particle size distribution. Hot gas, from the rotary kiln as well as the calciner, heats the cement raw meal in the pre-heating cyclones. Before entering the rotary kiln, the raw meal enters the calciner. Calcination starts at about 600$^\circ$C in some of the pre-heating cyclones, while the main calcination occurs in the calciner that operates at 900-1100$^\circ$C. Typically, 90\% of the limestone is calcined when leaving the calciner and entering the rotary kiln.

\begin{figure}[tb]
    \centering
    \begin{subfigure}[b]{0.24\textwidth}
    \centering
        \begin{tikzpicture}        
            \draw[pattern=north east lines, pattern color=brown!50] (2,2) circle (2 cm);
            \draw[pattern=north east lines, pattern color=orange!50] (2,2) circle (1.8 cm);
            \draw[thick,-] (2,2) circle (1.5 cm);
            \fill[white] (2,2) circle (1.5 cm);
            \fill[blue!40!white] (2,2) circle (0.1 cm);
    
            \draw[thick,<->] (2,2) -- node[anchor=south]{$r_c$} (3.05,3.05) ;
            \draw[thick,<->] (2,2) -- node[anchor=south]{$r_r$} (3.8,2) ;
            \draw[thick,<->] (2,2) -- node[anchor=south]{$r_w$} (0.6,0.6) ;

            \draw[-] (3.8,2.5) -- (4,3.0) node[anchor=south]{$A_w$};
            \draw[-] (3.5,1.5) -- (4.5,1.5) node[anchor=south]{$A_r$};
            \draw[-] (2.0,3.0) -- (3.5,4.0) node[anchor=west]{$A_c$};
    
            \draw[thick,->] (1.1,2.5) -- node[anchor=south west]{$Q_{cr}$} (0.5,2.8) ;
            \draw[thick,->] (2,0.4) --  (2,0.05) node[anchor=north]{$Q_{rw}$};
            \draw[thick,->] (3.2,0.5) -- node[anchor=north west]{$Q_{we}$} (3.6,0.5) ;
        \end{tikzpicture}    
        \caption{Cross section}
        \label{fig:cross}
    \end{subfigure}%
    \begin{subfigure}[b]{0.24\textwidth}
    \centering
    \resizebox{1\textwidth}{!}{%
    \begin{tikzpicture}   
        \draw[pattern=north east lines, pattern color=brown!50] (0,1.0) -- (0,3) -- (0.5,4) -- (0.6,4) -- (0.1,3) -- (0.1,1) -- (0.6,0) -- (0.5,0)-- cycle;
        \draw[pattern=north east lines, pattern color=orange!50] (0.1,1.0) -- (0.1,3) -- (0.6,4) -- (0.7,4) -- (0.2,3) -- (0.2,1) -- (0.7,0) -- (0.6,0)--cycle;
        
        \draw[pattern=north east lines, pattern color=brown!50] (2.2,1.0) -- (2.2,3) -- (1.7,4) -- (1.8,4) -- (2.3,3) -- (2.3,1) -- (1.8,0) -- (1.7,0)-- cycle;
        \draw[pattern=north east lines, pattern color=orange!50] (2.1,1.0) -- (2.1,3) -- (1.6,4) -- (1.7,4) -- (2.2,3) -- (2.2,1) -- (1.7,0) -- (1.6,0)--cycle;

        \draw[thick,-] (2.25,3) -- node[anchor=south west,scale=0.5]{Wall} (2.5,3.0) ;
        \draw[thick,-] (2.15,2.5) -- node[anchor=south west,scale =0.5]{Refractory} (2.5,2.5) ;

        \draw[thick,->] (-0.5,-0.1) -- (-0.5,4.5)node[anchor= west,scale =0.5]{y};
        \draw[thick,-] (-0.6,-0.0) -- (-0.4,0)node[anchor= west,scale =0.5]{0};
        \draw[thick,-] (-0.6,4) -- (-0.4,4)node[anchor= west,scale =0.5]{$h_{tot}$};
        \draw[thick,-] (-0.6,3) -- (-0.4,3)node[anchor= west,scale =0.5]{$h_{cu}$};
        \draw[thick,-] (-0.6,1) -- (-0.4,1)node[anchor= west,scale =0.5]{$h_{cl}$};

        \draw[thick,->] (1.2,3.5) -- (1.2,4.2)node[anchor= west,scale =0.5]{gas mixture};

        \draw[thick,->] (0.7,-0.1) node[anchor= east,scale =0.5]{Fuel}  -- (0.9,-0.1) -- (0.9,0.2);
        \draw[thick,->] (1.1,-0.3) node[anchor= north,scale =0.5]{Kiln gas}  -- (1.1,0.2);
        \draw[thick,->] (1.6,-0.1) node[anchor= west,scale =0.5]{3rd Air} -- (1.5,-0.1) -- (1.5,0.2);
        \draw[thick,->] (1.6,-0.3) node[anchor= west,scale =0.5]{Cyclone flow} -- (1.3,-0.3) -- (1.3,0.2);

        \draw[-] (-0.6,2.4)node[anchor= east,scale =0.5]{$y_{k+\frac{1}{2}}$} -- (2.3,2.4);
        \draw[-] (-0.6,2.0)node[anchor= east,scale =0.5]{$y_{k-\frac{1}{2}}$} -- (2.3,2.0);

        \draw[->] (1.25,1.9) -- node[anchor=north east,scale =0.5]{$\Tilde{H}_{g,k-\frac{1}{2}}$} node[anchor=north west,scale =0.5]{$N_{g,k-\frac{1}{2}}$} (1.25,2.1);
        \draw[->] (1.25,2.3) -- node[anchor=south east,scale =0.5]{$N_{g,k+\frac{1}{2}}$}node[anchor=south west,scale =0.5]{$\Tilde{H}_{g,k+\frac{1}{2}}$} (1.25,2.5);

        \draw[->] (0.70,2.2)  arc (0:300:0.1cm) node[anchor=south west,scale=0.5]{$R_k$};

        \draw[<->] (1.22,0) -- node[anchor=north,scale =0.5]{$r_l$} (1.6,0);
        \draw[<->] (1.22,4) -- node[anchor=north,scale =0.5]{$r_u$} (1.6,4);
        
    \end{tikzpicture}
    }
    \caption{Axial profile, with segment \\notation for the $k$-th volume.}
    \label{fig:axial}
    \end{subfigure}
    \caption{Diagrams of the calciner profiles. The diagrams illustrates the dimensions and flow direction.}\label{fig:profiles}
\end{figure}
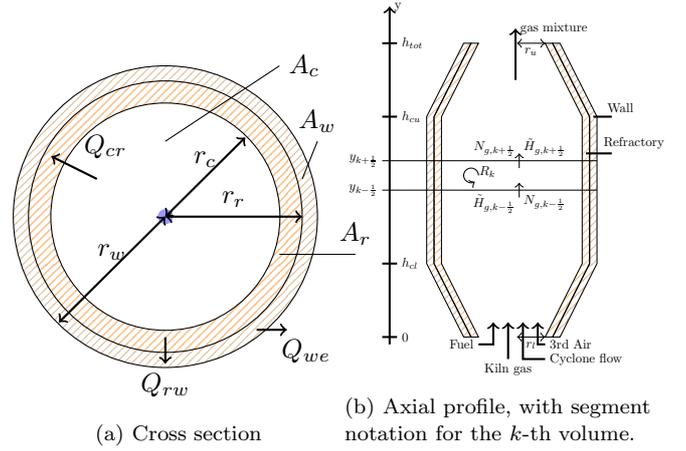

Fig. \ref{fig:profiles} illustrates the chamber of the calciner modeled in this paper. Fuel, gas, and cement raw meal enters at the bottom of the calciner and exit at the top. The fuel, the hot kiln gas, and the hot gas from the cooler provide the heat for the calcination. Different designs of varying complexity exist for calciners. 
The designs 
can be pipes, cyclones, or have another geometry. We do not model all of these types of calciners in this paper. The calciner modeled in this paper is a cylindrical chamber with a cylindrical cone in the top and the bottom. The calciner has a total height of $h_{tot}$ and is a cylinder with radius $r_c$ between the heights $h_{cl}$ and $h_{cu}$. The cone sections have the smaller radii, $r_l$ and $r_u$ for the lower and upper cone, respectively. 


\section{A mathematical model for dynamical simulation of the calciner}
\label{sec:CalcinerModel}

The calciner model is formulated as an index-1 DAE system.
The states, $x$, are the molar concentrations of each compound in the solid-gas mixture, $C$, and the internal energy densities of each phase, $\hat{U}$. The phase temperatures, $T$, and the pressure, $P$, are the algebraic variables, $y$. The resulting model can be represented as
\begin{subequations}\label{eq:dyn}
\begin{align}
    \partial_tx& = f(x,y;p),\quad x=[C;\hat U],\\
    0 &= g(x,y;p),\quad y=[T;P],
\end{align}
\end{subequations}
where $p$ is the system parameters. Manipulated variables and disturbances enter through the boundary conditions. The model is obtained using a systematic modeling methodology that integrates the a) geometry, b) thermo-physical properties, c) transport phenomena, d) stoichiometry and kinetics, e) mass and energy balances, and f) algebraic relations for the volume and the internal energy. The phases 
considered are the mixture (c), the refractory wall (r), and the shell wall (w). 


\subsection{Calciner and cement chemistry}
We use the standard cement chemist notation for the following compounds:
 \ce{(CaO)_2SiO_2} as \ce{C_2S}, \ce{(CaO)_3SiO_2} as \ce{C_3S}, \ce{(CaO)_3Al_2O_3} as \ce{C_3A}, and \ce{(CaO)_4} \ce{(Al_2O_3)(Fe_2O_3)} as \ce{C_4AF}, where C = \ce{CaO},   A = \ce{Al_2O_3}, S = \ce{SiO_2}, and  F = \ce{Fe_2O_3}. 

We use a finite-volume approach to describe the calciner in $n_v$ segments of length $\Delta y = h_{tot}/n_v$. We define the molar concentration vector, $C$, as mole per segment volume $V_\Delta(k)$, and assume all gasses are ideal. We employ the following assumptions: 1) The horizontal planes are assumed homogeneous with only dynamics along the height of the calciner (1D-model); 2) The gas flow is assumed to prevent solids from exiting the calciner through the bottom; 3) Only the 5 main clinker formation reactions are included; and 4) Only basic fuel reactions are included.

\subsection{Geometry}
The volume of segment $k$ is 
\begin{align}
    V_{\Delta}(k) = \pi r_c^2h_c(k) &+ \frac{\pi}{3}(r_c^2+r_u^2(k) +r_cr_u(k))h_u(k)\nonumber\\ &+ \frac{\pi}{3}(r_c^2+r_l^2(k) +r_cr_l(k))h_l(k).
\end{align}
The cylinder and cone heights are given by
\begin{subequations}
\begin{align}
    h_c(k) &= (y_{k+\frac{1}{2}}-h_u(k))-(y_{k-\frac{1}{2}}+h_l(k))\\
    h_u(k) &= \max(y_{k+\frac{1}{2}}-h_{cu},0),\\
    h_l(k) &= \max(h_{cl}-y_{k-\frac{1}{2}},0).
\end{align}
\end{subequations}
The small radius' of the cone sections are given by
\begin{subequations}
\begin{align}
    r_u(k) &= r_u+\frac{r_c-r_u}{h_{cu}}h_u(k),\\
    r_l(k) &= r_l+\frac{r_c-r_l}{h_{cl}}h_l(k).
\end{align}
\end{subequations}
The segment volumes for the refractory and walls are computed by the relations,
\begin{subequations}
\begin{align}
    V_{r,\Delta}(k) &= V_{\Delta}(k)|_{r_c=r_r} - V_{\Delta}(k),\\
    V_{w,\Delta}(k) &= V_{\Delta}(k)|_{r_c=r_w} - V_{\Delta}(k)|_{r_c=r_r}.
\end{align}
\end{subequations}
Each $V_{\Delta}(k)$ is computed with the radius for segment $k$.

Similarly, the surface area (sides) of each segment depends on location and is
\begin{equation}
\begin{split}
    A_c(k) &= 2 \pi r_c h_c(k)\\
    & + \pi (r_c + r_u(k)) \sqrt{h_u^2(k) + (r_c-r_u(k))^2}\\
    & + \pi (r_c + r_l(k)) \sqrt{h_l^2(k) + (r_c-r_l(k))^2}.
\end{split}
\end{equation}

\subsection{Thermo-physical model}
We provide a thermo-physical model for the enthalpy, $H(T,P,n$), and the volume, $V(T,P,n)$, of each phase. These models are homogeneous of order 1 in the mole vector, $n$. The thermo-physical model for $H$ and $V$ is
\begin{small}
 \begin{align}
    H(T,P,n) &= \sum_i n_i \left(\Delta H_{f,i}(T_0,P_0) + \int^T_{T_0}c_{p,i}(\tau)d\tau \right),
    \\
    V(T,P,n) &= \begin{cases} \sum_i n_i \left( \frac{M_i}{\rho_i} \right), \quad \text{solid, (s)},\\ \sum_{i} n_i \left( \frac{RT}{P} \right), \quad \text{gas, (g)}.
  \end{cases}    
\end{align}
\end{small}    
$\Delta H_{f,i}(T_0,P_0)$ is the formation enthalpy at standard conditions $(T_0,P_0)$. $M$ is the molar mass. As $H$ and $V$ are homogeneous of order 1, the enthalpy and volume density can be computed as
\begin{subequations}
\begin{align}
    \hat{H}_s &= H_s(T_c,P,C_s),\quad  &\hat{V}_s &= V_s(T_c,P,C_s),\\
    \hat{H}_g &= H_g(T_c,P,C_g),\quad &\hat{V}_g &= V_g(T_c,P,C_g),\\
    \hat{H}_r &= H_r(T_r,P,C_r),\quad &\hat{H}_w &= H_w(T_w,P,C_w).
\end{align}
\end{subequations}
For a given section, the solid and gas volumes can be obtained from their densities,
\begin{align}
    V_g = \hat{V}_gV_{\Delta}(k), \quad  V_s = \hat{V}_sV_{\Delta}(k).
\end{align}

\subsection{Transport phenomena}
In the calciner, mass is transported by convection (advection) and diffusion, while energy is transported by convection, diffusion, and radiation.

\subsubsection{Velocity:}
We assume that all material move uniformly (same speed and direction) and that the velocity is below 0.2 Mach \citep{Darcy-Howel}. In this case, the velocity of the turbulent flow of the mixture, $v_c$, can be described by the Darcy-Weisbach equation,
\begin{align} 
    v_c &= \Big(\frac{2}{0.316}\sqrt[4]{\frac{D_H^{5}}{\mu_m\rho_m^3}}\frac{|\Delta P|}{\Delta z}\Big)^{\frac{4}{7}}\text{sgn}\Big(\frac{-\Delta P}{\Delta z}\Big).
\end{align}
$\mu_m$ is the viscosity of the mixture, $\rho_m$ is the density of the mixture, $D_H$ is the hydraulic diameter for a non-uniform and non-circular cross-section channel \citep{HESSELGREAVES20171}. $\rho_m$ and $D_H$ are computed by
\begin{align}    
    \rho_m = \sum_iM_iC_{i}, \quad D_H &= \frac{4V_\Delta}{A_{c}}.   
\end{align}

\subsubsection{Viscosity and conductivity:} 
For a pure component gas, a correlation for the temperature-dependent viscosity is \citep{Sutherland1893}
\begin{align}
    \mu_{g,i} = \mu_0 \left(\frac{T}{T_0}\right)^{\frac{3}{2}}\frac{T_0+S_\mu}{T+S_\mu}.
\end{align}
$S_\mu$ can be calibrated given two measures of viscosity as in Table \ref{tab:Data-Coeff-gas}.

For a gas mixture, \cite{Wilke1950} provides a viscosity correlation, $\mu_g$, and the Wassiljewa equation with the Mason-Saxena modification provides a conductivity correlation, $k_g$ \citep{Poling2001Book}
\begin{subequations}
\begin{align}
    \mu_g &= \sum_i\frac{x_i\mu_{g,i}}{\sum_jx_j\phi_{ij}},\quad k_g = \sum_i\frac{x_ik_{g,i}}{\sum_jx_j\phi_{ij}},\\
    \phi_{ij} &= \bigg(1+\sqrt{\frac{\mu_{g,i}}{\mu_{g,j}}}\sqrt[4]{\frac{M_j}{M_i}}\bigg)^2\bigg(2\sqrt{2}\sqrt{1+\frac{M_i}{M_j}}\bigg)^{-1}.
\end{align}
\end{subequations}
$x_i$ is the mole fraction of component $i$.
The viscosity of the suspended gas mixture, $\mu_m$, is the given by the extended Einstein equation of viscosity \citep{TODA2006} 
\begin{align}
    \mu_m &= \mu_g\frac{1+\phi/2}{1-2\phi}, \ \phi = \frac{V_s}{V_\Delta}=\hat{V}_s.  \label{eq:mue}
\end{align}
Assuming that the solid-gas mixture can be considered as layers,  the thermal conductivity of the solid-gas mixture, $k_m$, is given by the serial thermal conductivity  \citep{Perry}
\begin{align}
    \frac{1}{k_m} = \frac{V_{g}}{V_\Delta}\frac{1}{k_g}+ \sum_i\frac{V_{s,i}}{V_\Delta}\frac{1}{k_{s,i}}.
\end{align}
The volumetric ratios describe the layer thickness.

\subsubsection{Mass transport:} 
The mass transport in the vertical direction is by convection (advection) and diffusion. The material flux vector is  
\begin{align}
    N = N_{a} + N_{d}, \quad N_a = v C, \quad N_d = -D \odot \partial_z C.
\end{align}

\begin{remark}
Note that the diffusion (dispersion) is low and set to zero in the simulations in this paper.
\end{remark}

\subsubsection{Enthalpy and heat transport:}
The vertical transport of of enthalpy (internal energy and pressure work) is given by the enthalpy flux that can be computed by
\begin{align}
    \Tilde{H}_c = H_g(T_c,P,N_g) + H_s(T_c,P,N_s)
\end{align}
The heat conduction is given by Fourier's law 
\begin{align}
    \Tilde{Q}_c = -k_m\partial_z{T}_c
\end{align}
with $k_m$ being the thermal conductivity.

\subsubsection{Heat transfer between phases:}
The transfer of heat between phases are described by Newton's law of heat transfer
\begin{subequations}
\begin{align}
    Q_{cr}^{cv} &= A_{cr}\beta_{cr}(T_c-T_r),\\
    Q_{rw}^{cv} &= A_{rw}\beta_{rw}(T_r-T_w),\\
    Q_{we}^{cv} &= A_{we}\beta_{we}(T_w-T_e).
\end{align}
\end{subequations}
$A_{ij}$ is the in-between surface area, and $\beta_{ij}$ is the convection coefficient. The convection coefficients, $\beta$, are computed by the correlation
\begin{align}
    \beta = \frac{k}{d}Nu,
\end{align}
where d is the diameter of the cross section. For thermal heat transport across surfaces of different phases, the overall heat transfer coefficient, $A \beta$, is given by
\begin{align}
    A\beta = \left(\frac{1}{A_0\beta_0} + \sum_{i=1}^{n-1}\frac{dx_i}{k_iA_i}+ \frac{1}{A_n\beta_n}\right)^{-1}.
\end{align}
$A_i$ is the surface area, $k_i$ is the conductivity, and $dx_i$ is the width of phase $i$. For curved surfaces the width is given by $dx_i = \ln(\frac{r_{i+1}}{r_i})r_i$.
The Gnielinski correlation can be used as a generic formula for the Nusselt number, $Nu$, of turbulent flow in tubes \citep{Incropera}
\begin{subequations}
\begin{align}
    Nu &= \frac{\frac{f}{8}(Re_D-1000)Pr}{1+12.7(\frac{f}{8})^\frac{1}{2}(Pr^\frac{2}{3}-1)},\\ f &=(0.79 \ln (Re_D) - 1.64)^{-2},\\
    Re_D &= \frac{\rho_m v_c D_H}{\mu_m} \\
    Pr &= \frac{C_{p} \mu_m}{k_m}.
\end{align}
\end{subequations}
The heat capacity $C_{p}$ is given by
\begin{align}
    C_{p} = \sum_in_ic_{p,i}
\end{align}
with $c_{p,i}$ being specific molar heat capacities.

\subsubsection{Radiation:}
The transfer of heat due to radiation is given by
\begin{subequations}
\begin{align}
    Q_{cr}^{rad} &= \sigma A_{cr}(\epsilon_cT^4_c-\epsilon_rT^4_r),\\
    Q_{we}^{rad} &= \sigma A_{we}(\epsilon_wT^4_w-\epsilon_eT^4_e).
\end{align}
\end{subequations}
$\sigma$ is the Stefan-Boltzmann's constant and $\epsilon$ is the emissivity.  We assume axial radiation is negligible. Assuming that the calciner is made of the same material as a rotary kiln, then the emissivity of the wall and refractory is 0.85 and 0.8 \citep{Hanein2017}. The total emissivity of the solid-gas mixture is given by \citep{ALBERTI2018274}
\begin{align}
    \epsilon_c = \epsilon_s + \epsilon_g - \Delta\epsilon^s_g, \quad \Delta\epsilon^s_g = \epsilon_s\epsilon_g,
\end{align}
where $\Delta\epsilon^s_g$ is the overlap emissivity. 
The emissivity of the gas mixture, $\epsilon_g$, is computed using the WSGG model of 4 grey gases \citep{Johanson2011}.
Assuming that the raw meal has the same emissivity as the kiln bed surfaces, $\epsilon_s = 0.9$ is the solid emissivity \citep{Hanein2017}.

\subsection{Stoichiometry and kinetics}
The stoichiometric matrix, $\nu$, and the reaction rate vector, $r = r(T,P,C)$, provide the production rates, $R$:
\begin{align}
    R = \begin{bmatrix}
        R_s\\R_g
    \end{bmatrix}&= \nu^Tr.
\end{align}
$R_s$ is the production rate vector of the solids (\ce{CaCO3}, \ce{CaO}, \ce{SiO2}, \ce{AlO2}, \ce{FeO2}, \ce{C_2S}, \ce{C_3S}, \ce{C_3A}, \ce{C_4AF}, and \ce{C}) and $R_g$ is the production rate vector of the gasses (\ce{CO_2}, \ce{N_2} \ce{O_2}, \ce{Ar}, \ce{CO}, \ce{H_2}, and \ce{H_2O}).
The reactions in the solid-liquid phase related to clinker production are
\begin{subequations}
\begin{alignat}{5}
    \text{$\#1$: }& & \ce{CaCO3} &\rightarrow \ce{CO2} + \ce{CaO}, \quad & r_1,\\
    \text{$\#2$: }& & 2\ce{CaO} + \ce{SiO_2} &\rightarrow \ce{C_2S}, \quad & r_2,\\
    \text{$\#3$: }& &\ce{CaO} + \ce{C_2S}&\rightarrow \ce{C_3S}, \quad & r_3,\\
    \text{$\#4$: }& &3\ce{CaO} + \ce{Al_2O_3}&\rightarrow \ce{C_3A}, \quad & r_4,\\
    \text{$\#5$: }& &4\ce{CaO} + \ce{Al_2O_3} + \ce{Fe_2O_3}&\rightarrow \ce{C_4AF}, \quad & r_5, 
    \end{alignat}
\end{subequations}
while the combustion of fuel reactions related to heat generation are  
\begin{subequations}
    \begin{alignat}{5}
   \text{$\#6$: }& & 2\ce{CO}+\ce{O_2}&\rightarrow 2\ce{CO_2}, \quad & r_6, \\
   \text{$\#7$: }& & \ce{CO} + \ce{H_2O} &\rightarrow \ce{CO_2} + \ce{H_2}, \quad & r_7,\\
   \text{$\#8$: }& & 2\ce{H_2}+\ce{O_2}&\rightarrow 2\ce{H_2O}, \quad & r_8, \\
   \text{$\#9$: }& &2\ce{C} +\ce{ O_2} &\rightarrow 2\ce{CO}, \quad & r_9,\\
   \text{$\#10$: }& &\ce{C} + \ce{H_2O}&\rightarrow \ce{CO} + \ce{H_2}, \quad & r_{10}, \\
   \text{$\#11$: }& &\ce{ C} + \ce{CO_2}&\rightarrow 2\ce{CO}, \quad & r_{11}.
\end{alignat}
\end{subequations}
The rate functions, $r_j(T,P,C)$, used in this paper are given given by expressions of the form
\begin{align}
    r = k(T)\prod_lP_l^{\beta_l}C_l^{\alpha_l},\quad k_j(T) = k_{0}T^ne^{-\frac{E_{A}}{RT}}.
\end{align}
 $k(T)$ is the modified Arrhenius expression. $C_l$ is the concentration (mol/L). $\alpha_l$ is either the stoichiometric related or experimental-based power coefficient. $\beta_l$ is the power of the partial pressure $P_l = (C_l/\sum_j C_j ) P$. 

\subsection{Mass and energy balances}
The mass balances for the solid phase and the gas phase are
\begin{subequations}
\begin{align}
    \partial_t {C}_{s} &= -\partial_zN_{s} + R_{s},\\
    \partial_t {C}_{g} &= -\partial_zN_{g} + R_{g}.
\end{align}
\end{subequations}
The energy balances for the combined solid and gas phases, the refractory wall , and the wall are
\begin{subequations}
\label{eq:EnergyBalances}
\begin{align}
    \partial_t\hat{U}_c &= -\partial_z(\Tilde{H}_c+ \Tilde{Q}_c) - \frac{Q^{rad}_{cr} + Q^{cv}_{cr}}{V_{\Delta}}, \\
    \partial_t\hat{U}_r &= - \partial_z\Tilde{Q}_r + \frac{Q^{rad}_{cr} + Q^{cv}_{cr}}{V_{r}} - \frac{Q^{cv}_{rw}}{V_{r}},\\
    \partial_t\hat{U}_w &= - \partial_z\Tilde{Q}_w + \frac{Q^{cv}_{rw}}{V_{w}} - \frac{Q^{rad}_{we} + Q^{cv}_{we}}{V_{w}} .
\end{align}
\end{subequations}

\subsection{Algebraic equations for volume and internal energy}
The total specific volume of the gas and the solid is governed by the relation
\begin{align}
    V_g(T_c,P,C_g) + V_s(T_c,P,C_s) = \hat{V}_{g} + \hat{V}_{s} = \hat{V}_{\Delta} = 1.
\end{align}
The specific energies, $\hat U$, in \eqref{eq:EnergyBalances}  can be related to temperature, pressure and concentration by the thermo-physical relations
\begin{small}
\begin{subequations}
\label{eq:EnergyAlgebra}
\begin{align}
\begin{split}
    \hat{U}_c &= \hat{H}_s + \hat{H}_g - P\hat{V}_g
    \\&=H_s(T_c,P,C_s) + H_g(T_c,P,C_g) - P 
    V_g(T_c,P,C_g), 
\end{split}    
    \\
    \quad \hat{U}_r &= \hat{H}_r = H_r(T_r,P,C_r), 
    \\ \hat{U}_w &= \hat{H}_w = H_w(T_w,P,C_w). 
\end{align}
\end{subequations}
\end{small}
\section{Simulation Results}
\label{sec:SimulationResults}
We simulate the calciner during 60 min of operation to demonstrate the simulation model. We consider a 33 m high calciner with an inner radius of 3.08 m  and refractory and shell thicknesses of 0.21 m and 0.01 m, respectively. The refractory is made of alumina brick and the shell is iron. The operation is set to match a 234 ton/h clinker production with a consumption of 400 Kcal/kg clinker corresponding to 12 ton fuel per hour.
The solid inflow is 67.7 kg \ce{CaCO3}/s, 5.4 kg  \ce{CaO}/s, 7.2 kg \ce{SiO2}/s, 1.2 kg  \ce{Al2O3}/s, 2.3 kg  \ce{Fe2O3}/s, and 11.2 kg  \ce{C2S}/s at 850$^\circ$C.
The fuel and air inflow is 3.4 kg carbon/s, 1.5 kg \ce{N2}/s, and 0.5 kg \ce{O2}/s  at 60$^\circ$C.
In addition, the calciner receives a kiln gas inflow of 
6.7 kg \ce{CO2}/s, 19.6 kg \ce{N2}/s, 1.9 kg \ce{O2}/s, 0.4 kg \ce{Ar}/s, and 0.2 kg \ce{H2O}/s at 1100$^\circ$C, and a 3rd air intake of air at 950$^\circ$C with 0.1 kg \ce{CO2}/s, 25.7 kg \ce{N2}/s, 7.9 kg \ce{O2}/s, 0.5 kg \ce{Ar}/s, and  0.2 kg \ce{H2O}/s. The temperature of the ambient environment is 25$^\circ$C. Appendix \ref{app:PhysicalProperties} provides the parameters and physical properties.

\subsection{Manual model calibration }\label{sec:tuning}

Selected reaction rates are manually calibrated to qaulitatively fit the model to operational data.
The calcination reaction rate, $r_1$, was adjusted by a factor of 270 to give a 90.0\% conversion of \ce{CaCO3} to \ce{CaO}. The main combustion reaction rates, $r_6$ and $r_9$, were adjusted to obtain a suitable outflow temperature in the range of 850-900$^\circ$C. We obtained an outlet temperature of 887.5$^\circ$C by adjusting $r_6$ by a factor of $5\times 10^5$ and $r_9$ by a factor of  $60$.

  \begin{figure}
     \centering
     \includegraphics[width=0.45\textwidth,trim={0.75cm 6.7cm 1.5cm 2.6cm},clip]{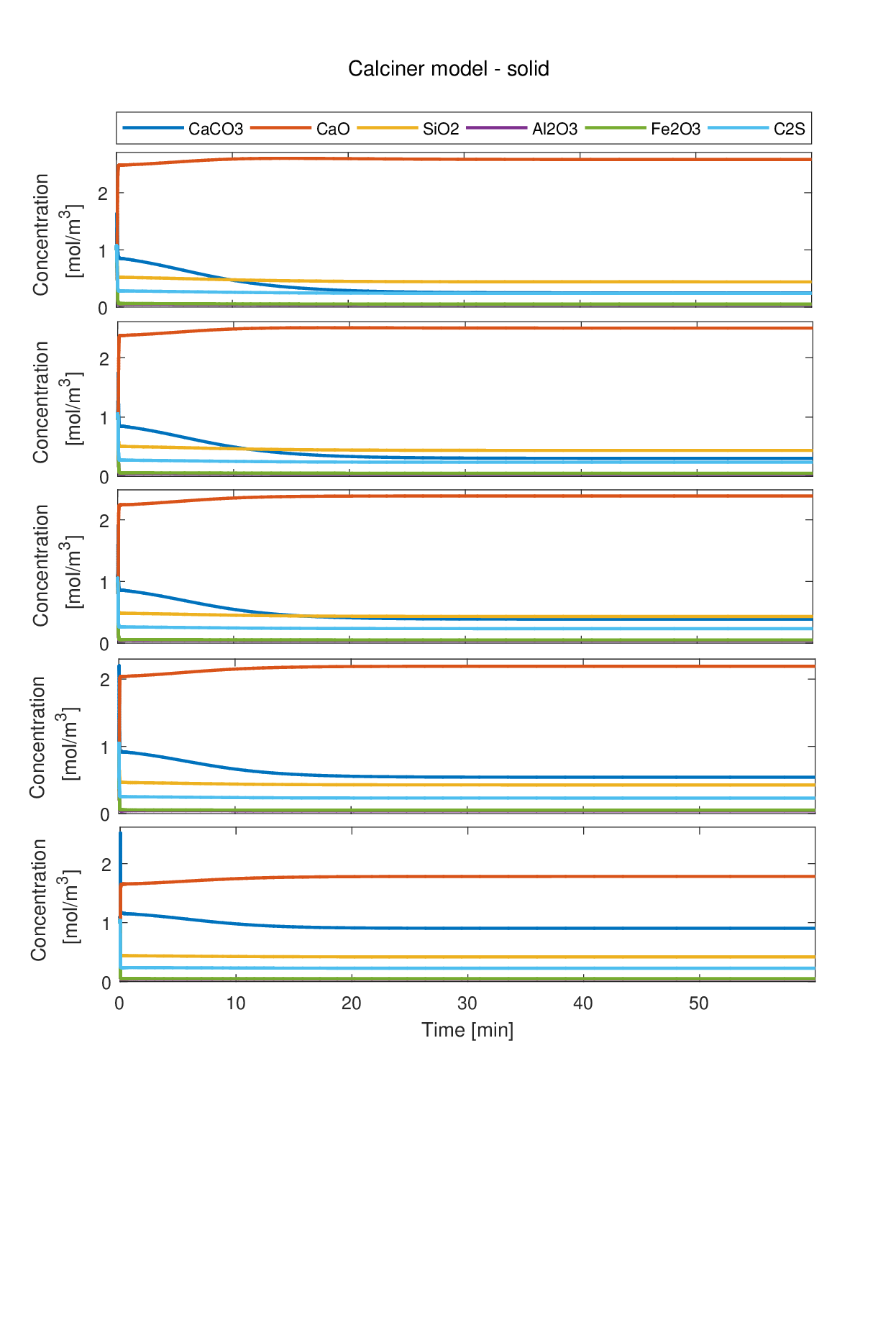}%
     \begin{subfigure}[b]{0.035\textwidth}
\begin{subfigure}[b]{\textwidth}
    \centering
    \resizebox{1\textwidth}{!}{%
    \begin{tikzpicture}   
       \draw[pattern=north east lines, pattern color=brown!50] (0,1.5) -- (0,4.5) -- (0.5,6) -- (0.6,6) -- (0.1,4.5) -- (0.1,1.5) -- (0.6,0) -- (0.5,0)-- cycle;
        \draw[pattern=north east lines, pattern color=orange!50] (0.1,1.5) -- (0.1,4.5) -- (0.6,6) -- (0.7,6) -- (0.2,4.5) -- (0.2,1.5) -- (0.7,0) -- (0.6,0)--cycle;
        
        \draw[pattern=north east lines, pattern color=brown!50] (2.2,1.5) -- (2.2,4.5) -- (1.7,6) -- (1.8,6) -- (2.3,4.5) -- (2.3,1.5) -- (1.8,0) -- (1.7,0)-- cycle;
        \draw[pattern=north east lines, pattern color=orange!50] (2.1,1.5) -- (2.1,4.5) -- (1.6,6) -- (1.7,6) -- (2.2,4.5) -- (2.2,1.5) -- (1.7,0) -- (1.6,0)--cycle;

        \draw[pattern=north east lines,pattern color=red] (0.3,4.8) -- (0.7,6) -- (1.6,6) -- (2.0,4.8) -- cycle;
        
        \node at (0.1,-0.4){};
    \end{tikzpicture}
    }
    \end{subfigure}  
    
    \begin{subfigure}[b]{\textwidth}
    \centering
    \resizebox{1\textwidth}{!}{%
    \begin{tikzpicture}   
        \draw[pattern=north east lines, pattern color=brown!50] (0,1.5) -- (0,4.5) -- (0.5,6) -- (0.6,6) -- (0.1,4.5) -- (0.1,1.5) -- (0.6,0) -- (0.5,0)-- cycle;
        \draw[pattern=north east lines, pattern color=orange!50] (0.1,1.5) -- (0.1,4.5) -- (0.6,6) -- (0.7,6) -- (0.2,4.5) -- (0.2,1.5) -- (0.7,0) -- (0.6,0)--cycle;
        
        \draw[pattern=north east lines, pattern color=brown!50] (2.2,1.5) -- (2.2,4.5) -- (1.7,6) -- (1.8,6) -- (2.3,4.5) -- (2.3,1.5) -- (1.8,0) -- (1.7,0)-- cycle;
        \draw[pattern=north east lines, pattern color=orange!50] (2.1,1.5) -- (2.1,4.5) -- (1.6,6) -- (1.7,6) -- (2.2,4.5) -- (2.2,1.5) -- (1.7,0) -- (1.6,0)--cycle;

        \draw[pattern=north east lines,pattern color=red] (2.0,4.8) -- (2.1,4.5) -- (2.1,3.6) -- (0.2,3.6) -- (0.2,4.5) -- (0.3,4.8) --  cycle;
        
        \node at (0.1,-0.4){};
    \end{tikzpicture}
    }
    \end{subfigure}  
    
    \begin{subfigure}[b]{\textwidth}
    \centering
    \resizebox{1\textwidth}{!}{%
    \begin{tikzpicture}   
        \draw[pattern=north east lines, pattern color=brown!50] (0,1.5) -- (0,4.5) -- (0.5,6) -- (0.6,6) -- (0.1,4.5) -- (0.1,1.5) -- (0.6,0) -- (0.5,0)-- cycle;
        \draw[pattern=north east lines, pattern color=orange!50] (0.1,1.5) -- (0.1,4.5) -- (0.6,6) -- (0.7,6) -- (0.2,4.5) -- (0.2,1.5) -- (0.7,0) -- (0.6,0)--cycle;
        
        \draw[pattern=north east lines, pattern color=brown!50] (2.2,1.5) -- (2.2,4.5) -- (1.7,6) -- (1.8,6) -- (2.3,4.5) -- (2.3,1.5) -- (1.8,0) -- (1.7,0)-- cycle;
        \draw[pattern=north east lines, pattern color=orange!50] (2.1,1.5) -- (2.1,4.5) -- (1.6,6) -- (1.7,6) -- (2.2,4.5) -- (2.2,1.5) -- (1.7,0) -- (1.6,0)--cycle;

        \draw[pattern=north east lines,pattern color=red] (0.2,2.4) -- (2.1,2.4) -- (2.1,3.6) -- (0.2,3.6) -- cycle;
        
        \node at (0.1,-0.4){};
    \end{tikzpicture}
    }
    \end{subfigure}  
    
    \begin{subfigure}[b]{\textwidth}
    \centering
    \resizebox{1\textwidth}{!}{%
    \begin{tikzpicture}   
       \draw[pattern=north east lines, pattern color=brown!50] (0,1.5) -- (0,4.5) -- (0.5,6) -- (0.6,6) -- (0.1,4.5) -- (0.1,1.5) -- (0.6,0) -- (0.5,0)-- cycle;
        \draw[pattern=north east lines, pattern color=orange!50] (0.1,1.5) -- (0.1,4.5) -- (0.6,6) -- (0.7,6) -- (0.2,4.5) -- (0.2,1.5) -- (0.7,0) -- (0.6,0)--cycle;
        
        \draw[pattern=north east lines, pattern color=brown!50] (2.2,1.5) -- (2.2,4.5) -- (1.7,6) -- (1.8,6) -- (2.3,4.5) -- (2.3,1.5) -- (1.8,0) -- (1.7,0)-- cycle;
        \draw[pattern=north east lines, pattern color=orange!50] (2.1,1.5) -- (2.1,4.5) -- (1.6,6) -- (1.7,6) -- (2.2,4.5) -- (2.2,1.5) -- (1.7,0) -- (1.6,0)--cycle;

        \draw[pattern=north east lines,pattern color=red] (0.2,1.5) -- (0.3,1.2) -- (2.0,1.2) -- (2.1,1.5) -- (2.1,2.4) -- (0.2,2.4) -- cycle;

        \node at (0.1,-0.4){};
    \end{tikzpicture}
    }
    \end{subfigure} 
    
    \begin{subfigure}[b]{\textwidth}
    \centering
    \resizebox{1\textwidth}{!}{%
    \begin{tikzpicture}   
        \draw[pattern=north east lines, pattern color=brown!50] (0,1.5) -- (0,4.5) -- (0.5,6) -- (0.6,6) -- (0.1,4.5) -- (0.1,1.5) -- (0.6,0) -- (0.5,0)-- cycle;
        \draw[pattern=north east lines, pattern color=orange!50] (0.1,1.5) -- (0.1,4.5) -- (0.6,6) -- (0.7,6) -- (0.2,4.5) -- (0.2,1.5) -- (0.7,0) -- (0.6,0)--cycle;
        
        \draw[pattern=north east lines, pattern color=brown!50] (2.2,1.5) -- (2.2,4.5) -- (1.7,6) -- (1.8,6) -- (2.3,4.5) -- (2.3,1.5) -- (1.8,0) -- (1.7,0)-- cycle;
        \draw[pattern=north east lines, pattern color=orange!50] (2.1,1.5) -- (2.1,4.5) -- (1.6,6) -- (1.7,6) -- (2.2,4.5) -- (2.2,1.5) -- (1.7,0) -- (1.6,0)--cycle;

        \draw[pattern=north east lines,pattern color=red] (0.3,1.2) -- (0.7,0) -- (1.6,0) -- (2.0,1.2) -- cycle;
        
        \node at (0.1,-2.2){};
    \end{tikzpicture}
    }
    \end{subfigure}  
\end{subfigure} 
     
     \caption{Evolution of the solid concentration in the calciner. The calciner is divided into 5 finite volumes. The top plot is the finite volume in the top of the calciner, while the bottom plot is the finite volume in the bottom of the calciner.}
     \label{fig:c1}
 \end{figure}
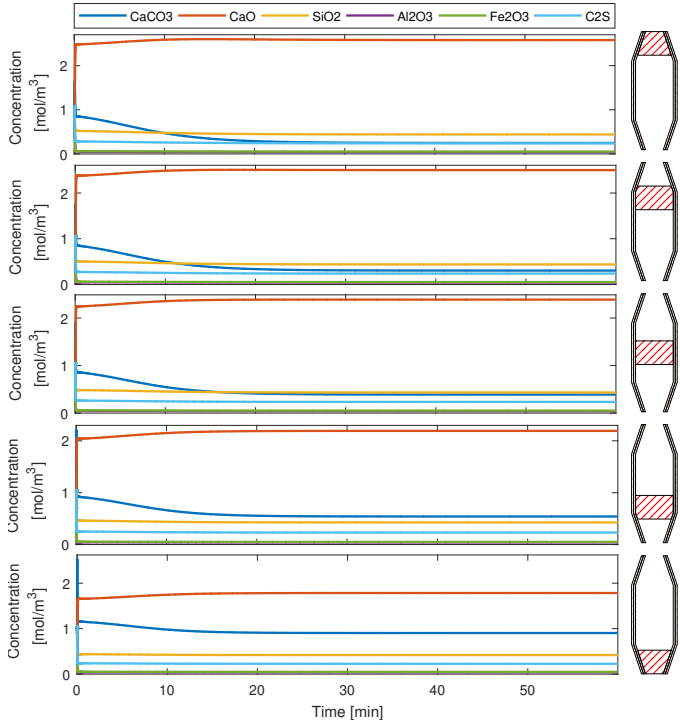
 
   \begin{figure}
     \centering
     \includegraphics[width=0.50\textwidth,trim={1.4cm 6.0cm 1.4cm 1.7cm},clip]{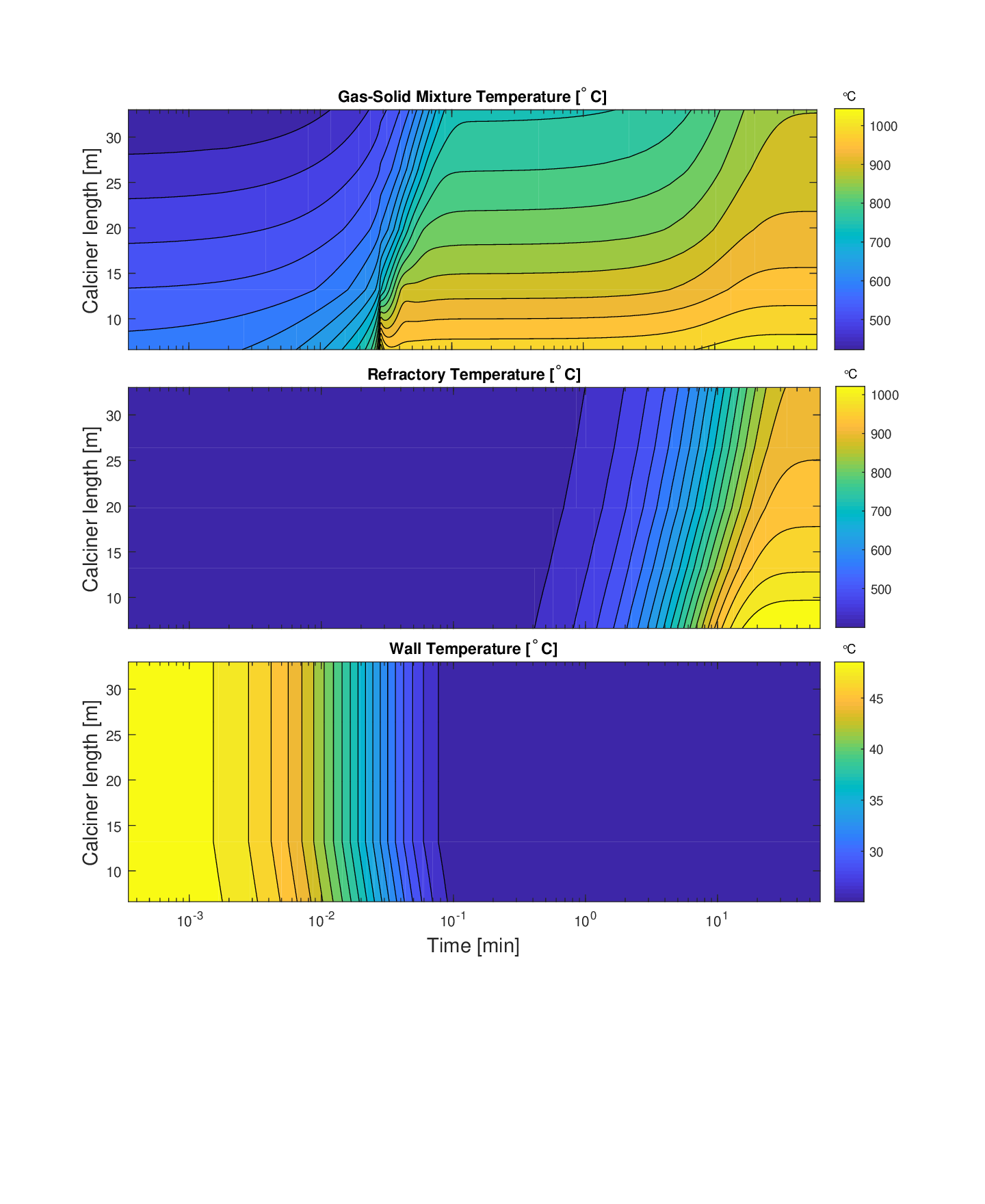}
     \caption{The temperature of the of the gas solid mixture, the refractory temperature, and wall temperature during the 60 min operation of the calciner. The temperatures are illustrated as function of position and time. Note that the time axis is logarithmic.}
     \label{fig:t1}
 \end{figure}
\subsection{Dynamic simulation}
Fig. \ref{fig:c1} and Fig. \ref{fig:t1} illustrate the dynamic behavior of the model. Fig. \ref{fig:c1} shows the solid molar concentrations and Fig. \ref{fig:t1} shows the temperatures. In this simulation, the system settles to a steady state in 30-40 min.
The mixture and refractory wall temperatures increase rapidly from 400$^\circ$C to 900-1000$^\circ$C along the calciner. Right after the ignition point at 1.6 s, the model exhibit non-monotonic behavior. Otherwise the behavior is monotonic approaching the steady states. The calcination process occurs when sufficient heat is released by the combustion.  

\subsection{Steady-state simulation}
\begin{figure}
     \centering
     \includegraphics[width=0.5\textwidth,trim={2.1cm 0.5cm 1.5cm 0.4cm},clip]{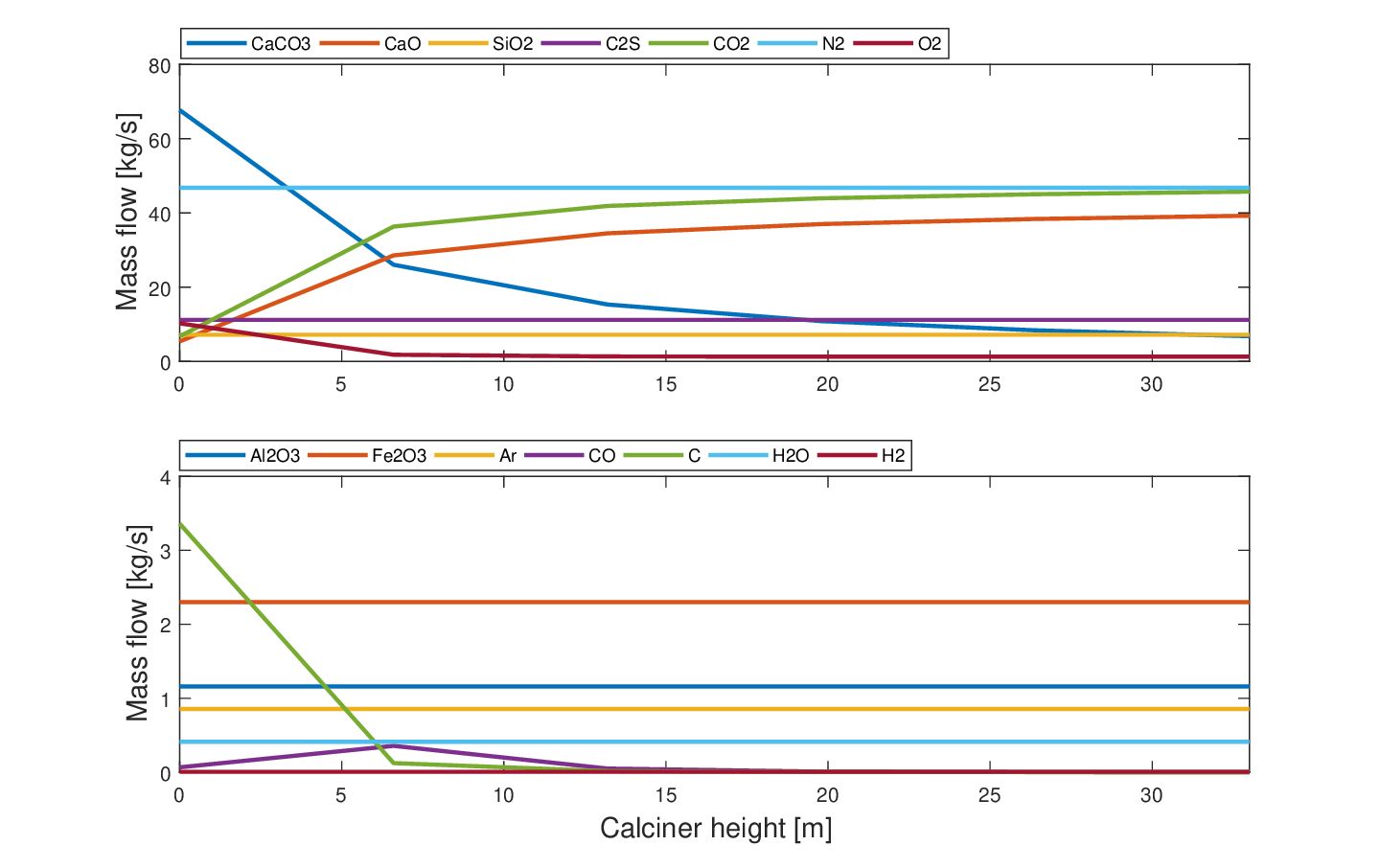}
     \caption{The steady-state mass flows of gasses and solids along the height of the calciner. Most of the calcination occurs in the first 10 m of the calciner.}
     \label{fig:n1}
 \end{figure}

Fig. \ref{fig:n1} shows the steady-state mass flow of all compounds in the calciner.
The calcination process occurs rapidly within the first 10 m where most of the \ce{CaO} is produced. Similarly, the combustion produces \ce{CO} from the fuel in the lower part of the calciner. This \ce{CO} is subsequently consumed in the \ce{CO2} forming reactions.
The \ce{CaCO3} conversion were obtained by manual calibration of the model parameters. The gas compositions can be used to evaluate the qualitative correctness of the steady states.
Typical practical operations conditions would be about $35\%$ \ce{CO2} and about $3\%$ \ce{O2} with the rest being primarily \ce{N2}. The model produced the steady-state outlet composition: 37.17\% \ce{CO2}, 59.74\% \ce{N2}, 1.42\% \ce{O2}, 0.77\% \ce{Ar}, 0.00\% \ce{CO}, 0.00\%\ce{C}, 0.82\% \ce{H2O}, and 0.08\% \ce{H2}. Hence the fuel is completely consumed and that the level of \ce{CO2} is in the right range (slightly higher than the practical operation). Correspondingly, the \ce{O2} level is also in the right range (slightly lower than the practical operation).

\section{Conclusion}\label{sec:Conclusion}
In this paper we presented an index-1 DAE model for the calciner in the pyro-section of a cement plant. The model can be used for dynamical and steady-state simulation of the calciner. By manual calibration, the model is able to provide simulation results that qualitatively matches practical operation. The model is the result of a systematic modeling procedure that involves flow pattern, geometry, thermo-physical properties, transport-phenomena, stoichiometry and kinetics, mass and energy balances, and algebraic relations for the volume and the internal energy. The calciner model will be connected with models for the pre-heating cyclones, the rotary kiln, and the cooler such that the pyro-section of a cement plant can be dynamically simulated. Such a model is important for model-based design and development of control and optimization systems for the pyro-section of cement plants.

\bibliography{ifacconf,biblio}             
                                                   







\newpage
\appendix
\section{Physical properties}\label{app:PhysicalProperties}

Table \ref{tab:reaction2}-\ref{tab:heatCap} provide the parameters and physical properties. Table \ref{tab:reaction2} shows the parameters for the kinetic expressions.
Table \ref{tab:Data-Coeff-solid} and Table \ref{tab:Data-Coeff-gas} shows literature data  for the solid and gas material properties. 

The molar heat capacity of \ce{CaCO3} is \citep{Jacobs1981}
\begin{equation}
\begin{split}
c_p &= -184.79 + 0.32 \cdot 10^{-3}T -0.13\cdot 10^{-5}T^2 \\ & -3.69\cdot 10^{6}T^{-2} + 3883.5T^{-\frac{1}{2}} \qquad  [\frac{\text{J}}{\text{mol}\cdot \text{K}}]
\end{split}
\end{equation}
for 298-750 K. 

The specific heat capacities of the remaining components are computed by \citep{Svensen2024Kiln}
\begin{equation}
c_p = C_0 + C_1 T + C_2 T^2
\end{equation}
 and Table \ref{tab:heatCap} reports the corresponding parameters ($C_0$, $C_1$,$C_2$).


\begin{table}[hb]
    \caption{Reaction rate coefficients of gasses from the open literature.}\label{tab:reaction2}%
    \begin{tabular*}{0.5\textwidth}{c| c| c|c | c | c |c | c | c }
    & Unit & $k_r $ & n  & $E_{A}$ & $\alpha_1$ & $\alpha_2$ & $\alpha_3$ & $\beta_2$\\ \hline
    $r_1$$^a$ &$\frac{\text{kg}}{\text{m}^3s}$&$10^{8}$ & 0 & $175.7$& $1$  & $0$  &0&0 \\
    $r_2$$^a$&$\frac{\text{kg}}{\text{m}^3s}$&$10^{7}$ & 0 & $240$& $2$  & $1$  &0 &0\\
    $r_3$$^a$&$\frac{\text{kg}}{\text{m}^3s}$&$10^{9}$ & 0 & $420$& $1$ & $1$ &0 &0\\
    $r_4$$^a$&$\frac{\text{kg}}{\text{m}^3s}$&$10^{8}$ & 0 & $310$& $3$ & $1$  &0&0 \\
    $r_5$$^a$&$\frac{\text{kg}}{\text{m}^3s}$&$10^{8}$ & 0 & $330$& $4$ & $1$ & $1$  &0 \\
        $r_6$$^b$&$\frac{\text{kg}}{\text{m}^3s}$&$7.0\cdot10^4$ & 0 & $66.5$& $1$\footnotemark[1]  & $1$\footnotemark[1]  &0 &0\\
        $r_7$$^c$ &$\frac{\text{mol}}{\text{m}^3s}$& $2.8\cdot10^6$ & 0 & 83.7 & 1& 1& 0&0\\
        $r_8$$^d$&$\frac{\text{mol}}{\text{m}^3s}$& $1.4\cdot10^{6}$ & 0.5 & $295.5$ & 1 & 1 & 0&0\\
        $r_9$$^e$& $\frac{\text{mol}}{\text{m}^3s}$& $8.8\cdot10^{11}$ & 0 &239 & $0.5$\footnotemark[1]  & $0.5$\footnotemark[1]  & 0 &0\\        
        $r_{10}$$^f$& $\frac{1}{\text{s}}$& $2.6\cdot10^8$ & 0 & 237 & 0 & 0 &0&  0.6\\
        $r_{11}$$^f$& $\frac{1}{\text{s}}$& $3.1\cdot10^6$ & 0 & 215 & 0 & 0  &0& 0.4\\
     \hline
    \end{tabular*}
    \footnotesize{All $\beta_1,\beta_3$ is zero and the unit of the activation energy $E_{A}$ is $[\frac{\text{kJ}}{\text{mol}}]$, $^1$ Unclear in source, $^a$ \cite{Mastorakos1999CFDPF}, $^b$ \cite{Guo2003}, $^c$ \cite{JONES1988}, $^d$ \cite{Karkach1999}, $^e$\cite{Walker1985}, $^f$\cite{BASU2018211}}
\end{table}

\begin{table}[hb]
    \centering
    \caption{Material properties of the solids.}
    \begin{tabular}{c|c|c|c}
    \hline
           &\shortstack{Thermal\\ Conductivity} & Density & \shortstack{Molar \\mass}\\ \hline
         Units    & $\frac{\text{W}}{\text{K}\cdot\text{m}}$ & $\frac{\text{g}}{\text{cm}^3}$ & $\frac{\text{g}}{\text{mol}}$ \\ \hline
         \ce{CaCO_3} &  2.248$^a$& 2.71$^b$  &100.09$^b$\\ 
         \ce{CaO}     & 30.1$^c$ &  3.34$^b$ &56.08$^b$\\ 
         \ce{SiO_2}   &  1.4$^{a,c}$& 2.65$^b$  &60.09$^b$\\ 
         \ce{Al_2O_3} &  12-38.5$^c$ 36$^a$& 3.99$^b$  &101.96$^b$\\ 
         \ce{Fe_2O_3} &  0.3-0.37$^c$& 5.25$^b$  &159.69$^b$\\ \hline
         \ce{C2S}     &  3.45$\pm$0.2$^d$& 3.31$^d$  &$172.24^g$\\ 
         \ce{C3S}     & 3.35$\pm$0.3$^d$ & 3.13$^d$ & 228.32$^b$\\ 
         \ce{C3A}     &  3.74$\pm$0.2$^e$& 3.04$^b$ & 270.19$^b$\\ 
         \ce{C4AF}    &  3.17$\pm$0.2$^e$& 3.7-3.9$^f$  &$485.97^g$ \\ \hline        
    \end{tabular}   
    
    \footnotesize{$^a$ from \cite{Perry}, $^b$ from \cite{CRC2022},\\ $^c$ from \cite{Ichim2018}, $^d$ from \cite{PhysRevApplied}, $^e$ from \cite{Du2021}, $^f$ from \cite{Portland},\\ $^g$ Computed from the above results} 
    \label{tab:Data-Coeff-solid}
\end{table}
\begin{table}[tb]
    \centering
    \caption{Material properties of the gasses.}
    \def\arraystretch{1.5}
    \begin{tabular}{c|c|c|c|c}    
    \hline
           &\shortstack{Thermal\\ Conductivity$^a$} & \shortstack{Molar\\ mass$^a$} & Viscosity$^a$ & \shortstack{diffusion\\ Volume$^b$}\\ \hline
        Units   & $\frac{10^{-3}\text{W}}{\text{K}\cdot\text{m}}$ & $\frac{\text{g}}{\text{mol}}$& $\mu \text{Pa}\cdot\text{s}$ & cm$^3$\\ \hline
         \ce{CO_2}    &\shortstack{\strut 16.77 (300K)\\ 70.78 (1000K) }  & 44.01  & \shortstack{\strut 15.0 (T=300K)\\ 41.18 (1000K) }& 16.3\\\hline
         \ce{N_2} & \shortstack{\strut 25.97(300K)\\  65.36(1000K) }  &28.014 &  \shortstack{ \strut 17.89(300K)\\  41.54(1000K) } & 18.5\\\hline
         
         \ce{O_2}  & \shortstack{\strut 26.49(300K)\\  71.55(1000K) } &  31.998 &  \shortstack{\strut 20.65 (300K)\\ 49.12 (1000K) }&  16.3\\\hline
         
         \ce{Ar}  & \shortstack{\strut 17.84 (300K)\\ 43.58 (1000K) }& 39.948&  \shortstack{\strut  22.74(300K)\\  55.69(1000K) }&  16.2\\\hline
         
         \ce{CO}  & \shortstack{\strut 25(300K)\\  43.2(600K) } & 28.010& \shortstack{\strut 17.8(300K)\\  29.1(1000K) } &  18\\\hline
         
         \ce{C_{sus}} & - & 12.011 & - &  15.9\\\hline
         
         \ce{H_2O} & \shortstack{\strut 609.50(300K)\\  95.877(1000K) }&18.015 &  \shortstack{\strut 853.74(300K)\\  37.615(1000K) }&   13.1\\\hline
         
         \ce{H_2} & \shortstack{\strut 193.1 (300K)\\ 459.7 (1000K) }&2.016 &  \shortstack{ \strut 8.938(300K)\\ 20.73 (1000K) }& 6.12\\\hline
    \end{tabular}
    $^a$ from \cite{CRC2022}, $^b$ from \cite{Poling2001Book}
    \label{tab:Data-Coeff-gas}
\end{table}

\begin{table}[hb]
    \centering
    \caption{Molar heat capacities.}
    \begin{tabular}{c|c |c | c| c}
            & $C_0$ & $C_1$ & $C_2$ & \shortstack{Temperature\\ range}\\ \hline
         Units & $\frac{\text{J}}{\text{mol}\cdot \text{K}}$& $\frac{10^{-3}\text{J}}{\text{mol}\cdot \text{K}^2}$&$ \frac{10^{-5}\text{J}}{\text{mol}\cdot \text{K}^3}$ & K\\\hline        
         \ce{CaO}$^b$     &  71.69& -3.08  & 0.22   & 200 - 1800\\ 
         \ce{SiO_2}$^b$   & 58.91 &  5.02 & 0 &844 - 1800\\ 
         \ce{Al_2O_3}$^b$ &  233.004&-19.59   &0.94   & 200 - 1800\\ 
         \ce{Fe_2O_3}$^c$ & 103.9  & 0 & 0  &-\\ \hline
         \ce{C2S}$^b$     &  199.6& 0  &0   &1650 - 1800\\ 
         \ce{C3S}$^b$     &  333.92&  -2.33&  0 &200 - 1800\\ 
         \ce{C3A}$^b$    & 260.58  & 9.58/2 & 0  &298 - 1800\\ 
         \ce{C4AF}$^b$  &  374.43& 36.4 & 0  & 298 - 1863\\ \hline
         \ce{CO_2}$^a$    & 25.98 &43.61 &-1.49  & 298 - 1500\\
         \ce{N_2}$^a$ & 27.31&5.19 &-1.55e-04  &298 - 1500\\ 
         \ce{O_2}$^a$ & 25.82&12.63 &-0.36   &298 - 1100\\ 
         \ce{Ar}$^a$ & 20.79 & 0 & 0   &298 - 1500\\ 
         \ce{CO}$^a$ & 26.87& 6.94  & -0.08   & 298 - 1500\\ 
         \ce{C_{sus}}$^a$& -0.45& 35.53 & -1.31  &298 - 1500\\ 
         \ce{H_2O}$^a$&30.89 & 7.86  &0.25  &298 - 1300\\ 
         \ce{H_2}$^a$&   28.95& -0.58&  0.19 & 298 - 1500\\ \hline 
    \end{tabular}
    
    \footnotesize{ $^a$ coefficient from  \cite{Jacobs1981}, $^b$ coefficients from \cite{HANEIN2020106043}, $^c$ from \cite{CRC2022}  }
    \label{tab:heatCap}
\end{table}

\end{document}